\documentclass{amsart}
\usepackage{style}

\title[Hyperbolic Handlebody Complements in 3-Manifolds]{Hyperbolic Handlebody Complements in 3-Manifolds}

\author[C. Adams]{Colin Adams}
\address{Department of Mathematics, Williams College, Williamstown, MA 01267}
\email{cadams@williams.edu}

\author [F. Gomez-Paz]{Francisco Gomez-Paz}
\address{Department of Mathematics, MIT,  77 Massachusetts Avenue, Cambridge, MA 02139-4307} 
\email{pjgomez@mit.edu}

\author[J. Kang]{Jiachen Kang}
\address{Department of Mathematics, University of Michigan, 530 Church St, Ann Arbor, MI 48109}
\email{jiachenk@umich.edu}

\author[L. Krause]{Lukas Krause}
\address{Department of Mathematics, University of California,
970 Evans Hall, Berkeley, CA 94720-3840}
\email{lukrau2002@gmail.com}

\author[G. Li]{Gregory Li}
\address{Department of Mathematics, Harvard University, Cambridge, MA 02138}
\email{gregoryli@college.harvard.edu}

\author[C. Marple]{Chloe Marple}
\address{Department of Mathematics, Pomona College, Claremont, CA 91711}
\email{ckme2022@mymail.pomona.edu}

\author[Z.Tan]{Ziwei Tan}
\address{Department of Mathematics, Bryn Mawr College,906 New Gulph Rd, Bryn Mawr, PA 19010}
\email{ztan2@brynmawr.edu}

\begin{document}

\begin{abstract} Let $M_0$ be a compact and orientable 3-manifold. After capping off spherical boundaries with balls and removing any torus boundaries, we prove that the resulting manifold \( M \) contains handlebodies of arbitrary genus such that the closure of their complement is hyperbolic. We then extend the octahedral decomposition to obtain bounds on volume for some of these handlebody complements. 
\end{abstract}

 \maketitle

 \section{Introduction}



A compact 3-manifold $M$ is said to be hyperbolic if after capping off any spherical boundary components and removing any torus boundaries, it admits a hyperbolic metric, which is to say a metric of constant sectional curvature $-1$. If the manifold has boundary components of genus greater than 1, we say it is tg-hyperbolic if the higher genus boundaries appear as totally geodesic surfaces in the hyperbolic metric.  

In \cite{Myers}, Myers proved that any compact orientable 3-manifold $M$ contains a knot such that the knot exterior is tg-hyperbolic. 
There are advantages to this. For instance, since volumes of manifolds are well-ordered, this tells us that every compact orientable 3-manifold $M$ has a unique minimal hyperbolic volume for any knot complement in $M$. Similalrly, the manifold inherits other hyperbolic invariants.

In this paper, we extend Myers' result to show the following:
\begin{restatable}{thm}{mainthm}\label{mainthm}
    Given a compact orientable 3-manifold $M$ with or without boundary, and a positive integer $n \geq 2$, $M$ contains a genus $n$ handlebody such that the  closure of its complement is tg-hyperbolic.
\end{restatable}



This implies the following result:
\begin{restatable}{cor}{generalthm}\label{generalthm}
Let $M$ be a compact orientable 3-manifold. Then for any finite sequence of positive integers $(a_1, a_2, \dots, a_n)$, there is a choice of $a_1$ solid tori, $a_2$ genus 2 handlebodies, $\dots$, $a_n$ genus $n$ handlebodies in $M$, disjoint from one another, such that the closure of their complement is tg-hyperbolic.
\end{restatable}

Note that it then follows from Theorem \ref{mainthm} that for any compact orientable 3-manifold $M$, and any integer $n \geq 1$, we can associate a volume $v_n^{min}(M)$, which is the least volume of the closure of the hyperbolic complement of a genus $n$ handlebody in $M$.  Thus, we obtain a hyperbolic volume spectrum associated to $M$. Other hyperbolic invariants can also be applied, even though the original manifold may not be hyperbolic.


In order to prove the main theorem, we use a result of Thurston, also used by Myers. A 3-manifold is said to be simple if it contains no properly embedded essential spheres, disks, annuli, or tori. In \cite{Thurston}, Thurston proved that a compact orientable simple 3-manifold with boundary is tg-hyperbolic. So it is enough to prove that there are no essential surfaces of these types.

The general strategy in proving Theorem \ref{mainthm} is similar to that of Myers. He showed that every such 3-manifold $M$ has a special handle decomposition with four 1-handles attached to every 0-handle. Within each 0-handle, he placed a true lover's tangle and connected them through the 1-handles to obtain a knot with simple complement.  

In Section 2, we replace the true lover's tangle in a $0$-handle by a simple knotted graph tangle as in Figure \ref{fig:knottedgraph1}. 
This increases the genus of the handlebody that is a neighborhood of the resulting knotted graph by one. Repeating the process for additional 0-handles and using the fact the number of $0$-handles can be made arbitrarily large, this shows \( M \) contains handlebodies of arbitrary genus with hyperbolic complements.

To conclude that the closure of the handlebody complement is simple, we build it as unions of pairs of 3-manifolds with an incompressible gluing surface in their boundaries.  In Section 2, we introduce the gluing lemma and a series of sufficient conditions for each pair of 3-manifolds along with the gluing surface to be simple. 

To check the gluing surfaces are incompressible as part of the conditions of the gluing lemma, we calculate the fundamental group of the central part in our knotted graph in Section \ref{figures/fundamentalgp}. Then we first use the gluing lemma inside each 0-handle in Section \ref{atoroidal} and we apply the gluing lemma again between 0-handles and the remaining handles to show that the entire handlebody complement is simple. 

In Section \ref{volumebounds}, we extend the octahedral decomposition for links in $S^3$ from \cite{DThurston} and links in $F \times I$ from \cite{AdamsCalderonMayer} to  handlebody complements in $F \times I$ and other 3-manifolds. This allows us to obtain various upper bounds on volume in terms of crossing number for spatial graphs that are retracts of the handlebodies.

\section {Preliminaries}\label{Preliminaries}

We start by introducing the following definitions:

\begin{definition}
A 3-manifold $M$ is \textit{irreducible} if every sphere contained in it bounds a ball, and is \textit{boundary-irreducible} if every properly embedded disk cuts a ball from the manifold.
 An annulus in $M$ is \textit{essential} if it is properly embedded, incompressible, and not boundary-parallel. A torus in $M$ is \textit{essential} if it is incompressible and not boundary-parallel.
    A $3$-manifold is \textit{simple} if it is irreducible, boundary irreducible, and contains no essential tori or annuli.

A manifold $M$ is called \textit{sufficiently large} if $M$ contains a properly embedded incompressible surface other than a sphere or disk. 
A compact orientable irreducible boundary-irreducible sufficiently large 3-manifold is called a \textit{Haken manifold}.
\end{definition}


\begin{definition}
A compact 3-manifold $M$ is \textit{tg-hyperbolic} if after capping off spherical boundary components, the complement of the torus components of $\partial M$ has a complete Riemannian metric with finite volume and constant sectional curvature $-1$ with respect to which the nontorus components of $\partial M$ are totally geodesic.
\end{definition}




\begin{theorem} [Thurston \cite{Thurston}]
Every simple Haken manifold is hyperbolic. 
\end{theorem}

\begin{definition}\label{irr pair}
A 3-manifold pair $(M, F)$ consists of a 3-manifold $M$ and a 2-manifold $F$ in $\partial M$.  We define $(M, F)$ to be an \textit{irreducible 3-manifold pair} if $M$ is irreducible and $F$ is incompressible.
\end{definition}

We refer to Waldhausen's definition of handle decomposition on p. 61 of  \cite{Waldhausen}. A handle decomposition of a 3-manifold is a decomposition into four collections of  balls. The first is denoted $N^0$ and made up of what we call 0-handles.  The second is denoted $N^I$ and made up of what we call 1-handles, each identified with $I \times D$.  The third is denoted $N^{II}$ and made up of what we call 2-handles, each identified with $D \times I$. The fourth is denoted $N^{III}$ and made up of what we call 3-handles. We assume the following conditions: 

\begin{enumerate}
    \item  $N^{III} \cap (N^0 \cup N^I \cup N^{II}) = \partial N^{III}$
    \item For any 1-handle, we have:
    
    \begin{enumerate} 
    \item $I \times D \cap N^0 = \partial I \times D$
    \item $I \times D \cap N^{II} = I \times d$ where $d$ is a collection of arcs in $\partial D$ 
    \end{enumerate}

\item For any 2-handle, we have: 
    \begin{enumerate} 
    \item $D \times I  \cap  N^0 = e_1 \times I$
    \item $D \times I \cap N^{I} = e_2 \times I$ where $e_1$ and $e_2$ are collections of arcs in $\partial D$ such that $e_1 \cup e_2 = \partial D$ 
    \end{enumerate}

    \item For each component of 
    $N^I \cap N^{II}$, the induced product structures agree.
    
\end{enumerate}

 A {\bf special handle decomposition} of a compact 3-manifold $M$ is a handle decomposition of $M$ such that the further conditions below apply:

 \begin{enumerate}
\item The intersection of any handle with any other handle or with $\partial M$ is either empty or connected.
\item Each 0-handle meets exactly four 1-handles and six 2-handles.
\item Each 1-handle meets exactly two 0-handles and three 2-handles. 
\item Each pair of 2-handles either
    \begin{enumerate}
    \item meets no common 0-handle or 1-handle, or
    \item meets exactly one common 0-handle and no common 1-handle, or 
    \item meets exactly one common 1-handle and two common 0-handles.
    \end{enumerate}

\item The complement of any 0-handle in the union of the 0-handles and 1-handles is connected.
\item  The union of any 0-handle with all the 2-handles and 3-handles is a handlebody that meets $\partial M$ in a disjoint collection of disks.
 \end{enumerate}
 
 \medskip

 A handle of index $k$ is denoted by $h_i^k$.
We utilize this lemma from Myers:

\begin{lemma}\label{specialhandledecom} \cite{Myers} Every compact orientable 3-manifold $M$ has a special handle decomposition.
\end{lemma}

\begin{figure}[htbp]
    \centering
    \includegraphics[width=0.4\linewidth]{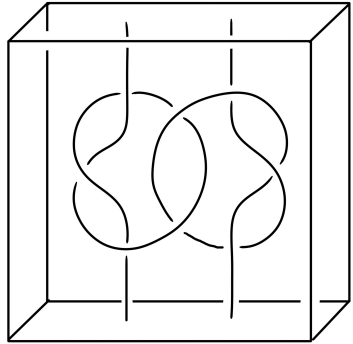}
    \caption{True Lovers Tangle $X$ with exterior $E(X)$ in a closed 3-ball}
    \label{fig:trueloverstangle}
\end{figure}

Having chosen a special handle decomposition, Myers properly embeds the tangle shown in Figure \ref{fig:trueloverstangle} into every 0-handle and connects the endpoints to the cores of the four 1-handles meeting each 0-handle. By choosing the connections appropriately the result is a knot. The true lover's tangle is simple in its tangle space \cite{Myers} and this allows Myers to prove that the knot complement is hyperbolic. 

Let $X$  denote the true lover's tangle and $E(X)$ denotes the exterior of $X$ in a closed 3-ball.

\begin{figure}[htbp]
    \centering
    \includegraphics[width=0.4\linewidth]{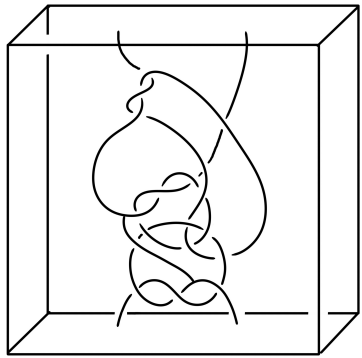}
    \caption{Knotted graph $J$ with exterior $E(J)$ in a closed 3-ball} 
    \label{fig:knottedgraph1}    
\end{figure}

Using a similar idea to Myers, we consider the knotted graph $J$ shown in Figure \ref{fig:knottedgraph1}. We will show that its exterior $E(J)$ in the closed ball is simple. 
This knotted graph is a 2-strand tangle with an arc attached.
 Thus, embedding this knotted graph into $(n-1)$ 0-handles and embedding a true lover's tangle into the remaining 0-handles, and then connecting the cores of the 1-handles results in an embedded graph in the manifold, a regular neighborhood of which is a genus $n$ handlebody. 

\begin{figure}[htbp]
    \centering
    \includegraphics[width=0.35\linewidth]{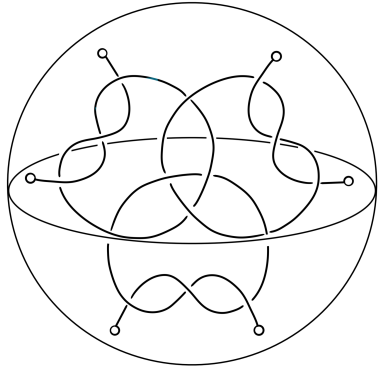}
    \caption{The knotted graph $I$ with the same exterior as $J$}
    \label{fig:knottedgraph2}
\end{figure}

\begin{figure}[htbp]
    \centering
    \includegraphics[width=1\linewidth]{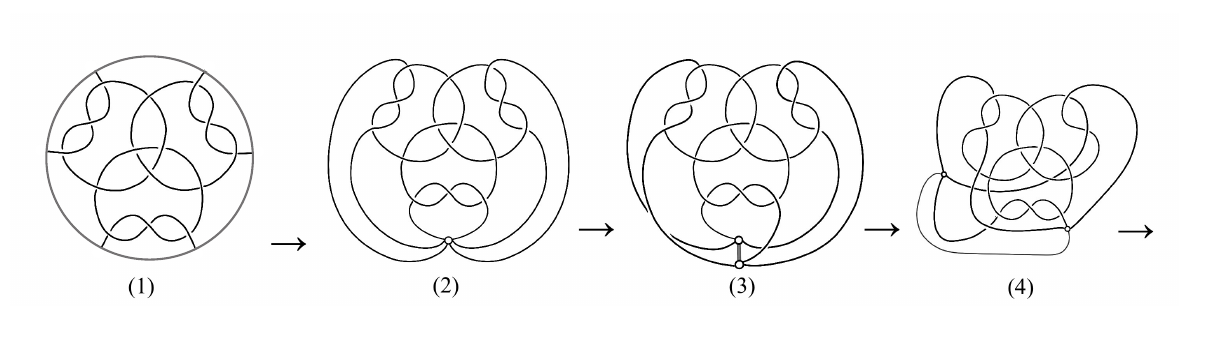}    \includegraphics[width=1\linewidth]{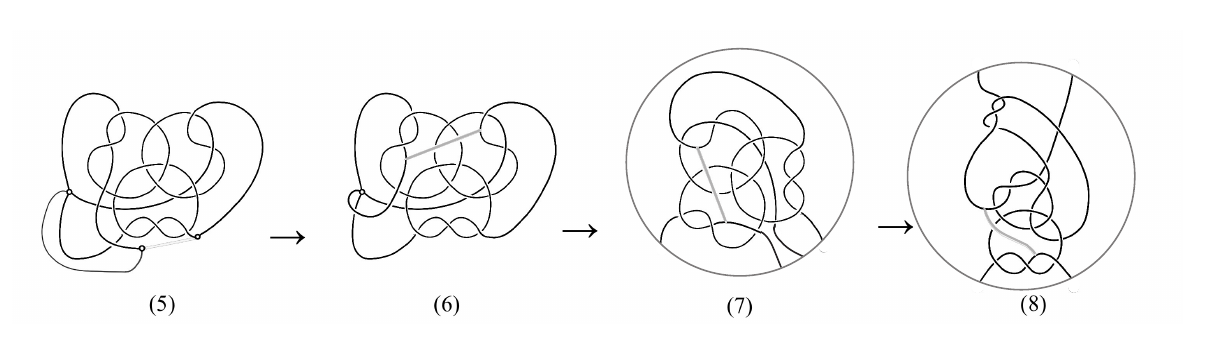}
    \caption{Showing the tangle exterior of $I$ and the spatial graph exterior of $J$ are equivalent}
    \label{fig:isotopy}
\end{figure}

\begin{remark}The knotted graph $J$ shown in Figure \ref{fig:knottedgraph1} is constructed from the three-strand tangle $I$ shown in Figure \ref{fig:knottedgraph2} thus their exteriors in a closed 3-ball are homeomorphic. Note that the 2-fold rotational quotient of the true lover's tangle gives the same tangle as the 3-fold rotational quotient of the tangle in Figure \ref{fig:knottedgraph2}, which is why we chose it as a candidate for being simple. 

The expression of $E(J)$ shown in Figure \ref{fig:knottedgraph2} is equivalent to the exterior of the graph formed by joining the six outgoing strands at a single vertex as appears in going from $(1)$ to $(2)$ in Figure \ref{fig:isotopy}. The subsequent moves in the figure  all preserve the exterior. 
Therefore, in subsequent sections we use the expression of $E(J)$ in Figure \ref{fig:knottedgraph2}.

\end{remark}

In Section \ref{atoroidal}, we prove that the exterior of this knotted graph in a closed 3-ball is hyperbolic using the following gluing gemma, and then use this to prove that the closure of the handlebody complement is simple.

 
Let $(M,F)$ be a pair of $3$-manifold $M$ and $2$-manifold $F\subset \partial M$. We use the following properties first introduced by Myers.

The idea of the gluing lemma is to decompose $M$ into two submanifolds that intersect exactly along a gluing face $F$, and prove each submanifold has certain properties that prevent essential spheres, disks, tori, and annuli in $M$.

\begin{definition}
    The pair $(M, F)$ has Property A if:
    \begin{enumerate}
        \item $(M, F)$ and $(M, \overline{\partial M - F})$ are irreducible $3$-manifold pairs
        \item No component of $F$ is a disk or $2$-sphere
        \item Every properly embedded disk $D$ in $M$ with $D \cap F$ a single arc is boundary-parallel.
    \end{enumerate}
    \medskip 

    The pair $(M,F)$ has Property B if:
    \begin{enumerate}
        \item $(M,F)$ has Property A
        \item No component of $F$ is an annulus
        \item Every properly embedded incompressible annulus $A$ in $M$ with $\partial A\cap \partial F=\varnothing$ is boundary-parallel.
    \end{enumerate}

    \medskip 

    The pair $(M,F)$ has property has Property C if:
    \begin{enumerate}
        \item $(M, F)$ has property B
        \item Every properly embedded disk $D$ in $M$ with $D \cap F$ a pair of disjoint arcs is boundary-parallel.
    \end{enumerate}

    \medskip

    The pair $(M,F)$ has Property B$'$ (respectively C$'$) if:
    \begin{enumerate}
        \item $(M, F)$ has Property B (respectively Property C)
        \item No component of $F$ is a torus
        \item Every incompressible torus in $M$ is boundary-parallel.
    \end{enumerate}
\end{definition}

Furthermore, we recall a useful gluing lemma proved by Myers in \cite{Myers}.

\begin{lemma}[Gluing Lemma]\label{lemma3.3-Myers} Let $M=M_0\cup M_1$, where each $M_i$ is a  compact orientable $3$-manifold and $F=M_0\cap M_1=\partial M_0\cap \partial M_1$ is a compact $2$-manifold. If $(M_0,F)$ has Property B$'$ and $(M_1,F)$ has Property C$'$, then $M$ is simple.
\end{lemma}

\begin{figure}[htbp]
    \centering
    \includegraphics[width=.6\linewidth]{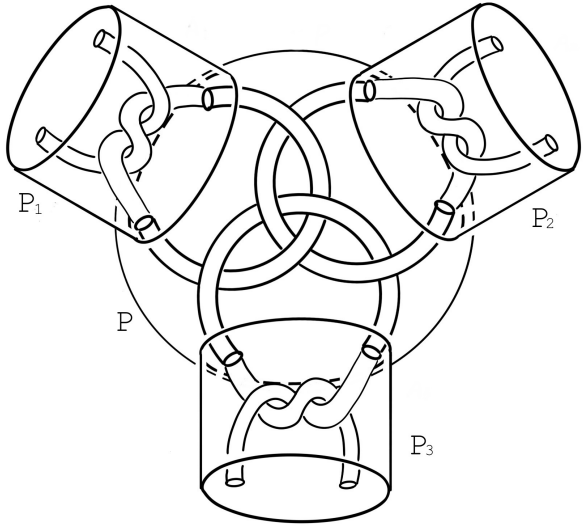}
    \caption{Decompose $E(I)$ into $P_1\cup P_2\cup P_3$ and $P$ with gluing faces $F_1,F_2,F_3$}
    \label{fig:decomp}
\end{figure}

In the following two sections, we prove $E(I)$ is simple using the gluing lemma. We first decompose it as in Figure \ref{fig:decomp}.

Let $F_i$ be the surface $P_i\cap P$. Then, $(P_i,F_i)$ has property B$'$ by Myers. So, by the gluing lemma, to prove $E(I)$ is simple it suffices to show that $(P,F_1\cup F_2\cup F_3)$ has property C$'$.
Now, $P$ is the exterior of the trivial $3$-tangle $T$ in the closed 3-ball, so we also denote $P$ by $E(T)$.

\begin{figure}[htbp]
    \centering
    \includegraphics[width=0.3\linewidth]{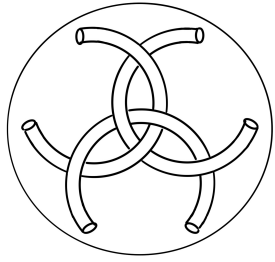}
    \caption{Trivial tangle $T$ with exterior $E(T)$ in the closed 3-ball $P$}
    \label{fig:tangle_exterior}
\end{figure}

We first show certain surfaces are incompressible in Section \ref{figures/fundamentalgp}, and then use these results to show $E(T)$ has Property C$'$ in Section \ref{atoroidal}.

\section {Fundamental Group Calculations}\label{figures/fundamentalgp}

\begin{figure}[htbp]
    \centering
    \includegraphics[width=0.4\linewidth]{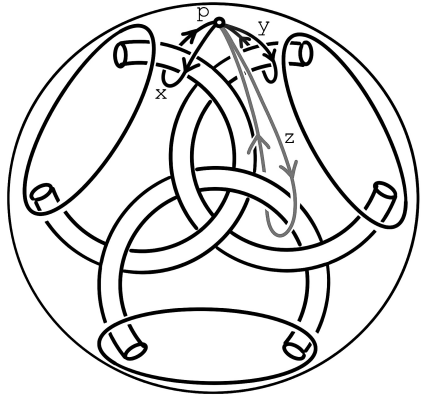}
    \caption{The three generators of the fundamental group $x,y,z$ of the exterior of our tangle}
    \label{fig:generator_diagram}
\end{figure}

We will show that various surfaces are $\pi_1$-injective in $E(T)$, which immediately implies incompressibility.  Pick a basepoint $p$ at the top of the ball in $E(T)$, and note that the fundamental group of this tangle exterior is free on three generators, $x, y$, and $z$. We've highlighted them in Figure \ref{fig:generator_diagram}.
To see that these are in fact generators, note that one can fix the endpoints of the handles they loop around and untangle $M$ into a handlebody while fixing these loops.


\begin{figure}[htbp]
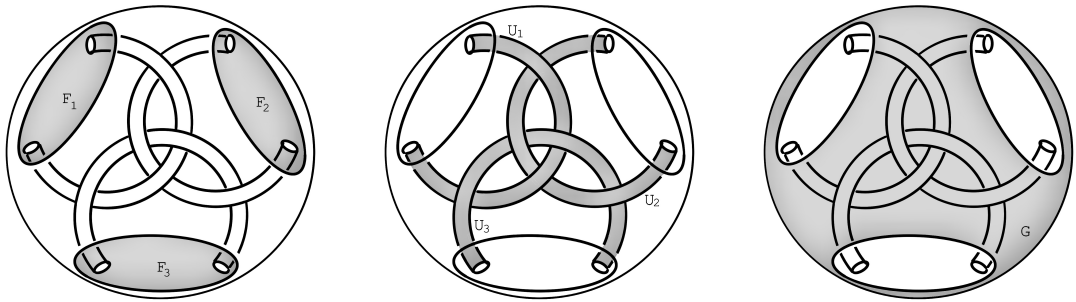

    \centering
    \includegraphics[width=0.3\linewidth]{figures/surface_F}
    \includegraphics[width=0.3\linewidth]{figures/Cmanifold_U}
    \includegraphics[width=0.3\linewidth]{figures/surface_G}
    \caption{Surfaces $F_1$, $F_2$, $F_3$, $U_1$, $U_2$, $U_3$ and $G$}
    \label{fig:surface_F}
\end{figure}

We highlight relevant surfaces $F_1, F_2, F_3$, $U_1, U_2, U_3$, and $G$ in Figure \ref{fig:surface_F}. 
Let $F := F_1 \cup F_2 \cup F_3$.

For any surface $S \subset E(T)$, let $i_*:\pi_1(S,p) \rightarrow \pi_1(E(T),p)$ denote the induced inclusion map on the fundamental groups.
We calculate the following maps in terms of $x, y$ and $z$.

\begin{enumerate}
    \item $i_*: \pi_1(F_1, p)  \rightarrow \pi_1(E(T),p)$
    \item $i_*: \pi_1(F_2,p)  \rightarrow \pi_1(E(T),p)$
    \item $i_*: \pi_1(G,p)  \rightarrow \pi_1(E(T),p)$
    \item $i_*: \pi_1(F_2 \cup U_2 \cup G, p) \rightarrow \pi_1(E(T),p)$ 
    \item $i_*: \pi_1(F_1 \cup F_2 \cup G, p)) \rightarrow \pi_1(E(T),p)$    
\end{enumerate}

\begin{remark}
     It is useful to pick our base point $p$ to lie at the top of the sphere.  However, in doing so $p$ will often not lie in the surfaces we wish to analyze. In these cases we must instead pick a basepoint $p'$ on the relevant surface $S$. After doing so the inclusion induces a map $i_*:\pi_1(S, p') \rightarrow \pi_1(E(T),p')$. Since $E(T)$ is path connected we may compose with some change of base point isomorphism $ \pi_1(E(T),p') \rightarrow \pi_1(E(T),p)$. This yields a map $\pi_1(S,p') \rightarrow \pi_1(E(T),p') \rightarrow \pi_1(E(T),p)$.
     By abuse of notation we refer to this map as $i_*: \pi_1(S,p) \rightarrow \pi_1 (E(T),p) $.
\end{remark}

\begin{claim} \label{sec:pi1_1}
    $i_*: \pi_1(F_1, p) \rightarrow \pi_1(E(T), p)$ is injective.
    \end{claim}

\begin{proof}
    Figure \ref{fig:genF1} depicts $i_*$ of the generators $e_1, e_2$ of $F_1$:

    \begin{figure}[htbp]
        \centering
        \includegraphics[width=0.4\linewidth]{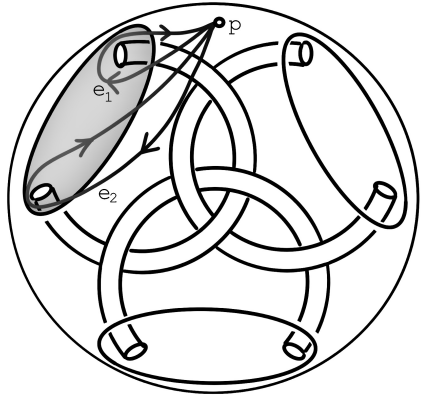}
        \caption{The generators $e_1$ and $e_2$ for $F$ appearing in $E(T)$}
        \label{fig:genF1}
    \end{figure}

    Note that $i_*(e_1)$ is easy and is just $x$. To see $i_*(e_2)$, we proceed by working backwards. 
    
    Starting with the loop $x^{-1}$, notice that this is the same as the loop in Figure \ref{fig:xinverse}.

    \begin{figure}[htbp]
        \centering
        \includegraphics[width=0.4\linewidth]{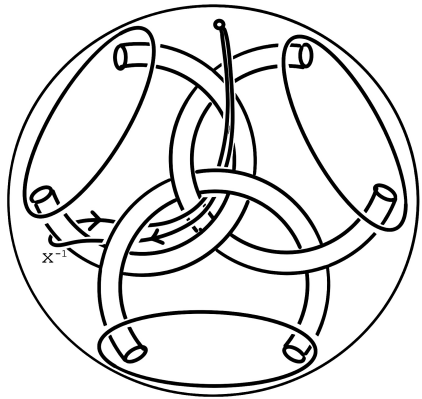}
        \caption{A realization of $x^{-1}$}
        \label{fig:xinverse}
    \end{figure}
    
    In order to obtain the loop $i_*(e_2)$, we need to change $x^{-1}$ so it passes over the strands of the tangle instead of under. This is achieved by conjugating. To go over the first strand, we conjugate by $z$ to get $zx^{-1}z^{-1}$, as depicted in Figure \ref{F1zxz}.

    \begin{figure}[htbp]
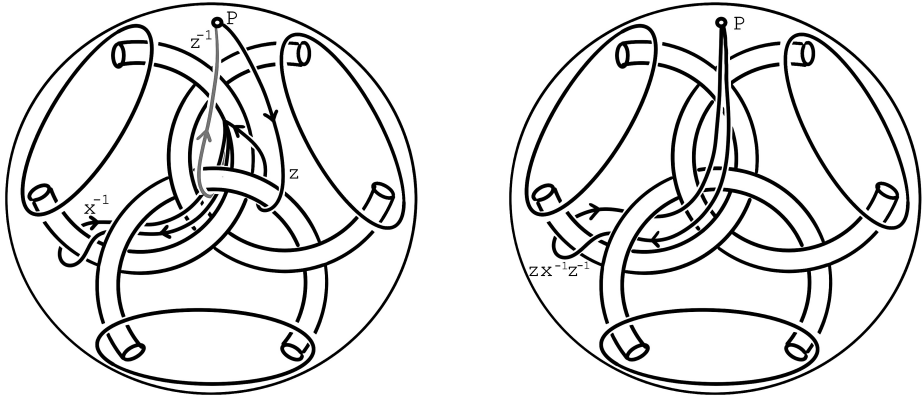

        \centering
        \begin{subfigure}{.4\textwidth}
            \centering
        \includegraphics[width=\linewidth]{figures/zxz1.pdf}
        \end{subfigure}
        \begin{subfigure}{.4\textwidth}
            \centering
           \includegraphics[width=\linewidth]{figures/zxz2.pdf}
        \end{subfigure}
        \caption{First step in constructing $i_*(e_2)$ in $i_*\pi_1(F_1,p)$}
        \label{F1zxz}
    \end{figure}

    To have it also go over the second strand, we first need to calculate the loop in Figure \ref{fig:step2}(A) in terms of the generators, then conjugate $zx^{-1}z^{-1}$ by it.

    \begin{figure}[htbp]
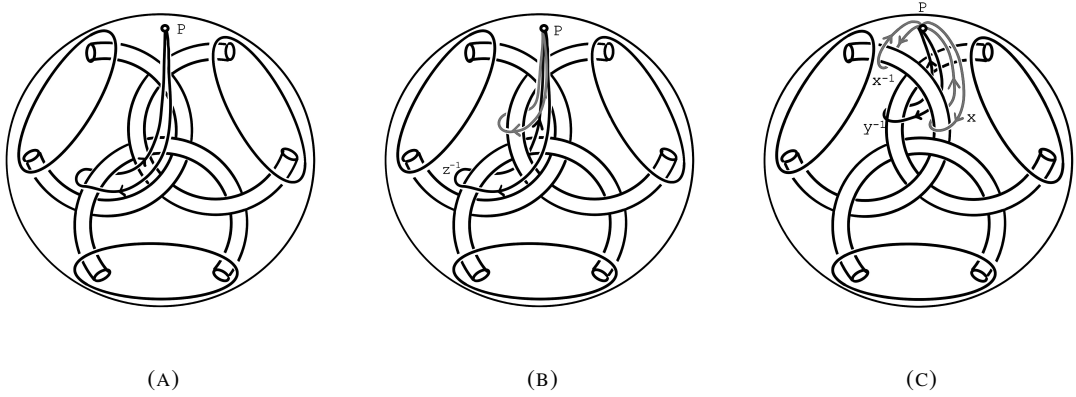

        \centering
        \begin{subfigure}{.3\textwidth}
            \centering
        \includegraphics[width=\linewidth]{figures/F1conj2.pdf}
        \caption{}
        \label{F1conj2}
        \end{subfigure}
        \begin{subfigure}{.3\textwidth}
            \centering
           \includegraphics[width=\linewidth]{figures/F1conj3.pdf}
           \caption{}
           \label{F1conj3}
        \end{subfigure}
        \begin{subfigure}{.3\textwidth}
            \centering
           \includegraphics[width=\linewidth]{figures/F1xyx.pdf}
           \caption{}
           \label{F1xyx}
        \end{subfigure}
        \caption{Second step in constructing $i_*(e_2)$ in $i_*\pi_1(F_1,p)$}
        \label{fig:step2}
    \end{figure}

    Notice that this loop is $z^{-1}$ conjugated by the gray loop in Figure \ref{fig:step2}(B), which is exactly $xy^{-1}x^{-1}$ as shown in Figure \ref{fig:step2}(C). This means the loop in Figure \ref{fig:step2}(A) is $(xy^{-1}x^{-1})z^{-1}(xyx^{-1})$.  Hence, we conjugate $zx^{-1}z^{-1}$ by this to obtain
    \[
    i_*(e_2) = [(xy^{-1}x^{-1})z^{-1}(xyx^{-1})]zx^{-1}z^{-1}[(xy^{-1}x^{-1})z(xyx^{-1})]
    \]
    
    Now, we claim that $i_*$ is injective. To see this, suppose $W$ is a non-trivial reduced word in $e_1$ and $e_2$. Notice that powers of $i_*(e_1) = x$ are reduced, and powers of $i_*(e_2)$ are conjugates of powers of $x$. Let $W$ be the word we're conjugating by; notice that binary products of $i_*(e_1)$ or its inverse with any words that are conjugated by $W$ or words conjugated by $W^{-1}$ are non-trivial. Therefore, the $i_*$-image of any non-trivial reduced word in $e_1$ and $e_2$ is still non-trivial in $x,y,z$, and thus the map $i_*$ is injective.
\end{proof}

\begin{claim} $i_*: \pi_1(F_2,p)  \rightarrow \pi_1(E(T),p)$ is injective.
\end{claim}

    \begin{proof}
    Although this is true by symmetry, since $E(T)$ has a rotational symmetry about a vertical axis in the plane sending $F_2$ to $F_1$, we will subsequently need calculations for the generators so we do that here.  We use the same method as in the proof of Claim \ref{sec:pi1_1}. Figure \ref{F2gens}(A) depicts $i_*(e_1)$ and $i_*(e_2)$, and we can see that $i_*(e_1)=y$. In Figure \ref{F2gens}(B) is the loop $y^{-1}$.

    \begin{figure}[htbp]
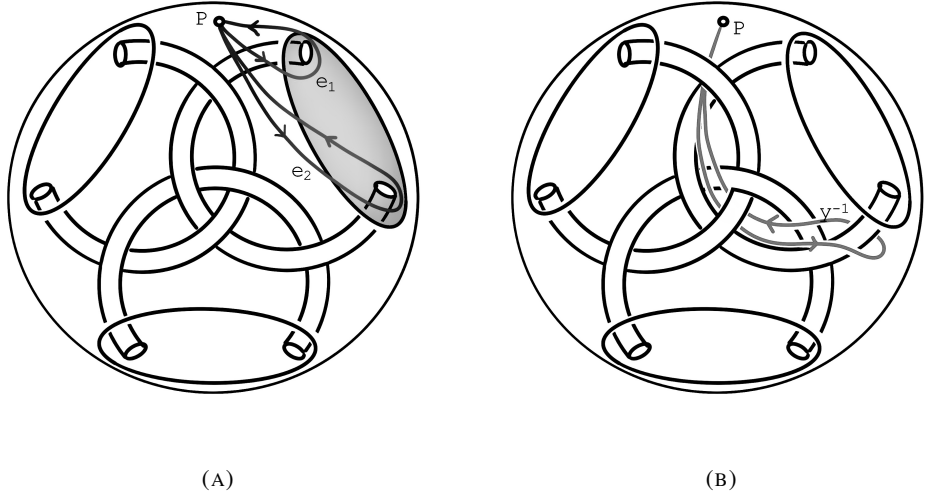

        \centering
        \begin{subfigure}{.4\textwidth}
            \centering
        \includegraphics[width=\linewidth]{figures/F2gen.pdf}
        \caption{}
        \end{subfigure}
        \begin{subfigure}{.4\textwidth}
            \centering
           \includegraphics[width=\linewidth]{figures/F2y1.pdf}
           \caption{}
        \end{subfigure}
        \caption{First step in constructing $i_*(e_2)$ in $i_*\pi_1(F_2,p)$}
        \label{F2gens}
    \end{figure}

    We first conjugate $y^{-1}$ by $x$ to go over the first crossing as shown in Figure \ref{F2zxzall}(A). Notice that the loop in Figure \ref{F2zxzall}(B) is $zx^{-1}z^{-1}$ as shown in Figure \ref{F2zxzall}(C), and we then conjugate $xy^{-1}x^{-1}$ by $zx^{-1}z^{-1}$ to go over the second crossing.

    \begin{figure}[htbp]
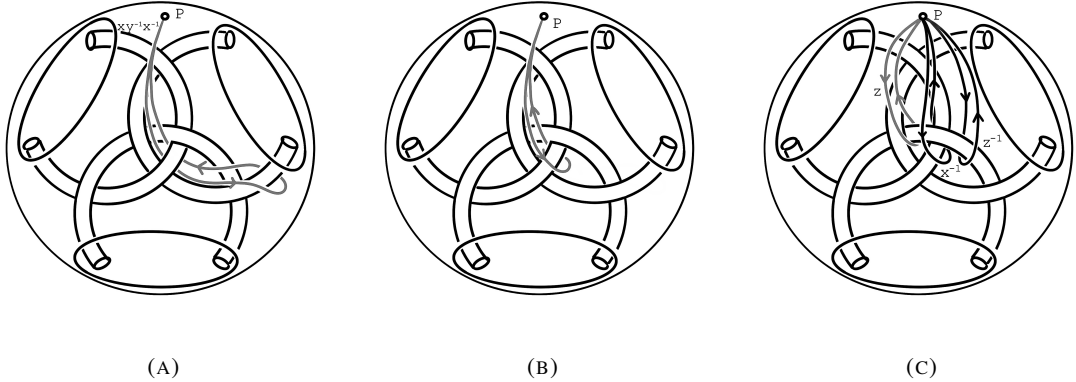

        \centering
        \begin{subfigure}{.3\textwidth}
            \centering
        \includegraphics[width=\linewidth]{figures/F2xyx.pdf}
        \caption{}
        \label{F2xyx}
        \end{subfigure}
        \begin{subfigure}{.3\textwidth}
            \centering
           \includegraphics[width=\linewidth]{figures/F2con2.pdf}
           \caption{}
           \label{F2conj2}
        \end{subfigure}
        \begin{subfigure}{.3\textwidth}
            \centering
           \includegraphics[width=\linewidth]{figures/F2zxz.pdf}
           \caption{}
           \label{F2zxz}
        \end{subfigure}
        \caption{Second step in constructing $i_*(e_2)$ in $i_*\pi_1(F_2,p)$}
        \label{F2zxzall}
    \end{figure}

    Thus, $i_*(e_1) = y$ and $i_*(e_2) = (zx^{-1}z^{-1})(xy^{-1}x^{-1})(zxz^{-1})$. Because the word $i_*(e_2)$ begins and ends with $z,z^{-1}$, any non-trivial reduced word in terms of $i_*(e_1)$ and $i_*(e_2)$ is still reduced in $x,y,z$.
    \end{proof}
    
\begin{claim} \label{Ginjects} $i_*: \pi_1(G,p))  \rightarrow \pi_1(E(T),p)$ is injective.
\end{claim}

\begin{figure}[htbp]
        \centering
        \begin{subfigure}{.4\textwidth}
         \centering
    \includegraphics[width=\linewidth]{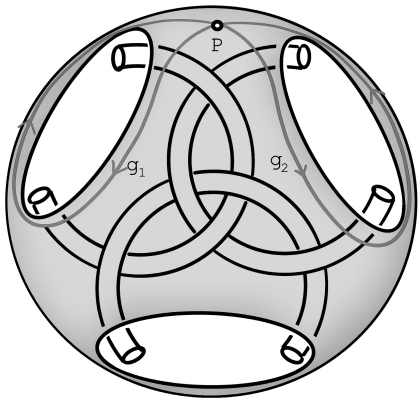}
    \caption{Generators of $G$.}
    \label{Ggen}
        \end{subfigure}
        \begin{subfigure}{.4\textwidth}
            \centering
    \includegraphics[width=\linewidth]{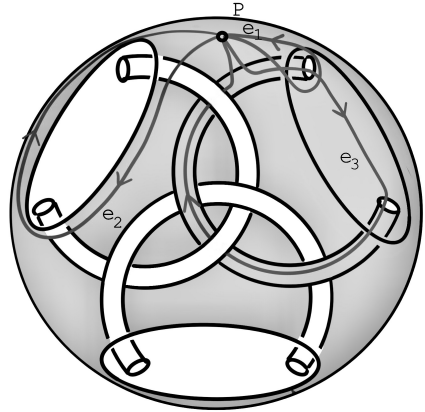}
    \caption{Generators of $F_2\cup U_2\cup G$.}
    \label{FUGgen}
        \end{subfigure}
        \caption{Generators for fundamental groups of surfaces}.
        \label{lotsofgenerators}
    \end{figure}

\begin{proof} Notice that $G$ is a thrice-punctured sphere. In particular, it has two generators.  If we choose the generators for $G$ to be $g_1$ and $g_2$ as shown in Figure \ref{lotsofgenerators}(A), then they are given by the product of the two generators for $F_1$ and the product of the two generators for $F_2$, namely $i_*(g_1) = xy^{-1}x^{-1}z^{-1}xyx^{-1}zx^{-1}z^{-1}xy^{-1}x^{-1}zxy$
and $i_*(g_2) = zx^{-1}z^{-1}xy^{-1}x^{-1}zxz^{-1}y$. Again, we note that non-trivial reduced words in $g_1$ and $g_2$ are also non-trivial in $x, y$, and $z$, so the map is injective, as desired.
\end{proof}

\begin{claim}\label{FUG}$i_*: \pi_1(F_2 \cup U_2 \cup G, p) \rightarrow \pi_1(E(T),p)$ is injective.
\end{claim}

\begin{proof}The surface $F_2 \cup U_2 \cup G$ is a twice-punctured torus 
with three generators, depicted in Figure \ref{lotsofgenerators}(B).
We've already calculated $i_*(e_1)=y$ and $i_*(e_2)=xy^{-1}x^{-1}z^{-1}xyx^{-1}zx^{-1}z^{-1}xy^{-1}x^{-1}zxy$ from Claim \ref{Ginjects}, so
it remains to calculate $i_*(e_3)$, which is the same as the loop in Figure \ref{FUGe}(A).


\begin{figure}[htbp]
        \centering
        \begin{subfigure}{.4\textwidth}
            \centering
        \includegraphics[width=\linewidth]{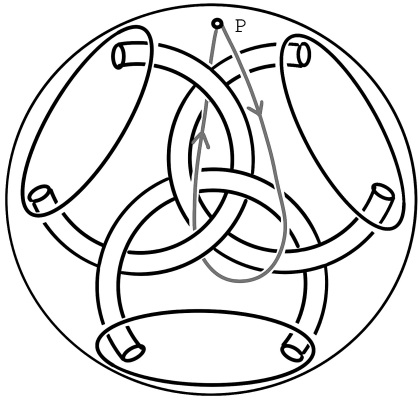}
        \caption{}
        \end{subfigure}
        \begin{subfigure}{.4\textwidth}
            \centering
           \includegraphics[width=\linewidth]{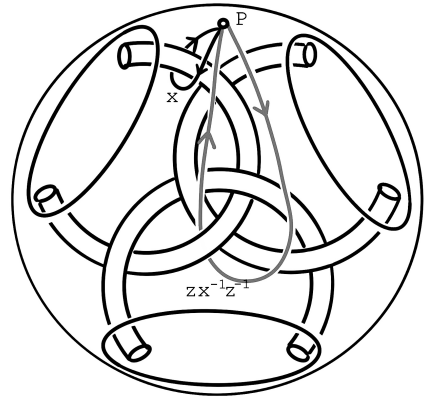}
           \caption{}
        \end{subfigure}
        \caption{Showing the fundamental group of $F_2\cup U_2\cup G$ injects}
        \label{FUGe}
    \end{figure}


    Because we have already shown $zx^{-1}z^{-1}$ is the loop from Figure \ref{F1zxz}, we can multiply this by $x$ as in Figure \ref{FUGe}(B) to obtain $i_*(e_3)=zx^{-1}z^{-1}x$.

Again, by inspection we notice that non-trivial words reduced in the $e_i$ are also non-trivial in their images, so $i_*$ is injective. By rotational symmetry of $E(T)$, this implies $F_1\cup U_1\cup G$ is also incompressible.
\end{proof}

\begin{claim}\label{FFG}$i_*: \pi_1(F_1 \cup F_2 \cup G, p)) \rightarrow \pi_1(E(T),p)$ is injective. 
\end{claim}

\begin{figure}[htbp]
    \centering
    \includegraphics[width=0.4\linewidth]{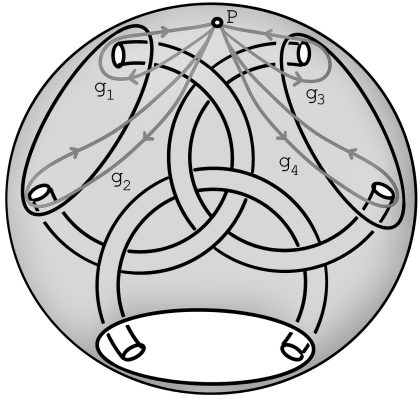}
    \caption{Generators for $F_1 \cup F_2 \cup G$}
    \label{F1F2G}
\end{figure}

\begin{proof}Notice that as in Figure \ref{F1F2G},  this is a five-punctured sphere, with generators $g_1,g_2,g_3,g_4$ given by the generators of $F_1$ and $F_2$. Thus, we have the following: $$i_*(g_1) = x,$$ $$i_*(g_2) =xy^{-1}x^{-1}z^{-1}xyx^{-1}zx^{-1}z^{-1}xy^{-1}x^{-1}zxyx^{-1},$$ 
$$i_*(g_3)=y,$$ $$i_*(g_4)=zx^{-1}z^{-1}xy^{-1}x^{-1}zxz^{-1}.$$

Because $\pi_1(F_1\cup F_2\cup G)=\langle g_i\mid 1\le i\le 4\rangle$, it's also generated by $g_1,g_2':=g_1g_3g_1^{-1}g_2g_1g_3^{-1}g_1^{-1},g_3,g_4$. Because there are exactly four of them, they are also a set of free generators.

Their images under $i_*$ are $x,z^{-1}xyx^{-1}zx^{-1}z^{-1}xy^{-1}x^{-1}z,y,zx^{-1}z^{-1}xy^{-1}x^{-1}zxz^{-1}$, so any reduced words in $g_1,g_2',g_3,g_4$ is reduced in $x,y,z$. This means $i_*$ is injective.
\end{proof}

\section{
Proof of Main Theorem}\label{atoroidal}


\begin{lemma}\label{ETFpA}
$(E(T), F)$ has Property $A$.
\end{lemma}

\begin{proof}
We verify that $(E(T),F)$ has the properties (1), (2), (3) needed for Property A.

\medskip

(1) $E(T)$ is a genus three handlebody. By Jaco \cite{Jaco} Example III.13(b), $E(T)$ has no incompressible $2$-spheres, i.e. is irreducible. The incompressibility of $F=F_1\cup F_2\cup F_3$ follows from \Cref{sec:pi1_1} and the rotational symmetry of $E(T)$, and incompressibility of $\overline{\partial M - F}=G\cup U_1\cup U_2\cup U_3$ follow from \Cref{Ginjects} and the fact each $U_i$ is obviously incompressible. 

\medskip

(2) All components of $F$ are disks with two punctures  and thus are neither disks nor $2$-spheres.

\medskip

(3) Let $D$ be a properly embedded disk in $E(T)$ whose boundary intersects $F$ in a single arc.
Notice the rotational symmetry of $E(T)$, so without loss of generality, we may assume $\partial D\cap F\subset F_1$.
Since $\partial D$ doesn't intersect $F_2$ nor $F_3$, $\partial D$ is also disjoint from $U_2$, $U_3$. Thus $\partial D\subset F_1\cup U_1\cup G$, which is incompressible by Claim \ref{FUG}. Thus there exists a disk $D' \subset F_1\cup U_1\cup G$ such that $\partial D' = \partial D$. Together $S = D \cup D'$ form a sphere. Recall any sphere in a handlebody bounds a ball. Thus $S$ bounds a ball, and hence we may isotope $D$ to $D'$ through the ball. Since $D' \subset \partial E(T)$, $D$ is boundary parallel as desired.
\end{proof}

\begin{lemma}\label{lem:annu_conj}
    Let $A_1,A_2$ be the two boundary components of a properly embedded annulus $A$ in a manifold $M$.
    Fix a base point $p\in M$ and choose paths $x$ and $z$ from $p$ to points on $A_1$ and $A_2$ respectively.  Conjugate each $A_i$ by its path to create two loops in $\pi_1(M,p)$.  Then these two loops are conjugate in $\pi_1(M,p)$.
\end{lemma}

\begin{proof}
Let $y$ be a path in $A$ from the endpoint of $x$ to the endpoint of $z$. Then $xA_1x^{-1}\simeq xyA_2y^{-1}x\simeq \\ (xyz^{-1})zA_2z^{-1}(zy^{-1}x^{-1})$, as we wished to show.
\end{proof}

\begin{lemma}
$(E(T), F)$ has Property B.
\end{lemma}

\begin{proof}

We verify that $(E(T),F)$ satisfies the conditions  (1), (2), and (3) for Property B.

\medskip

(1) This follows by Lemma \ref{ETFpA}.

\medskip

(2) All components of $F$ are twice punctured disks and therefore not annuli.

\medskip

(3) Fix some incompressible annulus $A$ such that $\partial A \cap \partial F = \varnothing$. We wish to show $A$ must be boundary parallel.
Denote the two boundary components of $A$ as $\partial A_1, \partial A_2$.
Note that since $\partial A \cap \partial F = \varnothing$, each $ \partial A_i$ must lie entirely in some $U_j, F_j$ or in $G$.

The only non-trivial simple closed curves on $U_j$ are the meridians of the strands. Thus if $\partial A_i$ lies in $U_j$ we may isotope it along $U_j$ so that it lies in $F_j$. Similarly, since $G$ is a thrice-punctured sphere, the only non-trivial simple closed curves in $G$ must wrap around some puncture. Thus if $\partial A_i$ lies in $G$ we may isotope it to lie in some $F_j$.

 Thus we need only consider the case where both $\partial A_1$ and $\partial A_2$ lie in $F$. We proceed by casework.

 \medskip

 \textbf{Case One:}{ $\partial A_1, \partial A_2$ lie in different components of $F$.}

 \medskip
 
 By rotational symmetry, we may without loss of generality assume $\partial A_1 \subset F_1$ and $\partial A_2 \subset F_2$.
 Note since $F_1$ is a twice-punctured disk, $\partial A_1$ must wrap around one puncture or around both punctures.
 
 Recall that one may compute the homology class of $\partial A_1$ by abelianizing its homotopy class in $\pi_1(E(T),p)$. 
 Thus the homology class of $\partial A_1$ is $x^{\pm 1}$ if it winds around one puncture and $0$ otherwise. Similarly, since $\partial A_2 \subset F_2$, the homology class of $\partial A_1$ is $y^{\pm 1}$ if it winds around one puncture and $0$ otherwise.

Note the annulus $A$ realizes a homology between $\partial A_1$ and $\partial A_2$. Thus they must have the same homology, which implies they must both wrap around both punctures. Thus the two loops can be isotoped to be the two generators of $\pi_1(G)$. Furthermore by Lemma \ref{lem:annu_conj}, the homotopy classes $[\partial A_1], [\partial A_2] \in \pi_1(E(T),p)$ must be conjugate. However these two elements are cyclically reduced and have different length, so by Theorem 1.3 of \cite{magnus}, this means they are not conjugate.

\medskip

\textbf{Case Two:}{ $\partial A_1, \partial A_2$ lie in the same component of $F$.}

\medskip

Without loss of generality,  we assume $\partial A_1, \partial A_2$ lie in $F_1$.
Since $\partial A_1$ and $\partial A_2$ are homologous they must both either wrap around one puncture, or both wrap around both punctures. Note if $\partial A_1$ and $\partial A_2$ wrap around different punctures, we may isotope $\partial A_1$ along $U_1$ to wrap around the same puncture as $\partial A_2$. Thus, without loss of generality, we can assume that $\partial A_1$ and $\partial A_2$ are non-trivial concentric simple closed curves on $F_1$.

Given $A$ is incompressible, then by Jaco \cite{Jaco} Example III.13 (b),  it must be boundary compressible. Let $D$ be the boundary compression disk with arcs $\alpha \subset A,\beta \subset \partial E(T)$. Since $\partial A_1$ and $\partial A_2$ are concentric simple closed curves on $F_1$, they bound an annulus $A'$ in $F_1$ and $\beta$ must  lie in $A'$. Then, surgering the torus given by $A \cup A'$ using the disk $D$ yields a sphere, which must bound a ball by irreducibility. Hence, $A \cup A'$ bounds a solid torus that can be used to push $A$ into $\partial E(T)$.

\end{proof}

\begin{lemma}
(E(T), F) has property C$'$.
\end{lemma}
\begin{proof}
We verify that $(E(T),F)$ has the properties (1), (2), (3), and (4) for Property C$'$.

\begin{enumerate}
    \item This follows directly from the lemmas above.
    
    \item All the components of $F$ have boundary and so are clearly not tori.
    
    \item By Jaco \cite{Jaco} Example III.13(b), the only incompressible and boundary-incompressible surface in a handlebody is a disk, so 
    $E(T)$ contains no essential tori.

    \item Let $D$ be a properly embedded disk in $E(T)$ whose boundary intersects $F$ in two disjoint arcs, $\alpha_1$ and $\alpha_2$. If the two arcs are both on $F_1$, then incompressibility of $F_1\cup U_1\cup G$ implies $D$ is boundary-parallel.

    If the two endpoints of $\alpha_1$ lie on the union of the the two hole boundaries in $F_1$, then to avoid  $\partial D$ being disconnected, $\alpha_2$ must also be on $F_1$. If one endpoint of $\alpha_1$ lies on a hole and the other endpoint lies on the outer boundary of $F_1$, for $\partial D$ to be a simple closed curve, it must intersect $\partial U_1$ in two points, which implies $\alpha_2$ also lies on $F_1$.

    So if $\alpha_i$'s are on distinct components of $F$, their endpoints must both lie on the outer boundary of $F$, which means they are disjoint from $U_i$'s. By incompressibility of $F_1\cup F_2\cup G$, $D$ is boundary-parallel.
\end{enumerate}
\end{proof}

\begin{lemma} \label{exteriorsimple} $E(I)$ is simple.
\end{lemma}

\begin{proof} Note that Myers proved that $(A_1 \cup A_2 \cup A_3, F)$ has property B$'$, and we have proven $(E(T), F)$ has property C$'$. Thus by Lemma \ref{lemma3.3-Myers}, $E(I)$ is simple.
\end{proof}


Since $E(J)$ is homeomorphic to $E(I)$, it is also simple. Let $\partial B$ be the 4-punctured sphere that bounds the ball containing $E(J)$ as in Figure \ref{fig:knottedgraph1}. We also prove the following lemma:

\begin{lemma}\label{ELpropertyC'}
    $(E(J),\partial B)$ has property C$'$.
\end{lemma}

\begin{proof}
    Because we have already proven that $E(J)$ is simple, it suffices to show that $\partial B$ and $\overline{\partial E(J)-\partial B}$ are incompressible in $E(J)$.


   We know $E(J)$ is boundary-irreducible.
    Both surfaces $\partial B$ and $\overline{\partial E(J)-\partial B}$ are spheres with four boundary components, each boundary component of which is nontrivial on $\partial E(J)$. 
   A nontrivial simple closed  curve in either $\partial B$ or $\overline{\partial E(J)-\partial B}$ is also nontrivial on $\partial E(J)$ because $\partial E(J)$ is homeomorphic to a genus $3$ surface.
   Since we have already proved that $E(J)$ is $\partial$-irreducible, such a curve is nontrivial in $E(J)$. So, both surfaces are incompressible.
   \end{proof}
   

 


We are now ready to prove the main Theorem:


\mainthm*
\begin{proof}

 If $\partial M$ contains spherical boundaries, we cap them off with balls.  Thus, without loss of generality, we suppose $\partial M$ has no sphere components. We further remove any torus boundaries.

    Let $C$ be a collar on $\partial M$ and let $N$ be $\overline{M-C}$. Using Lemma \ref{specialhandledecom}, we take a special handle decomposition of $N$.
    Note that Myers' proof of the lemma begins with a barycentric subdivision of a triangulation of the manifold, and then takes the dual cell complex. So, in the construction of the special handle decomposition, the number of $0$-handles is the same as $0$-cells in the dual cell complex. This means we may choose the original triangulation so that in the special handle decomposition of $N$, there are $m$ $0$-handles where $m > n$. 
 
    Let $H''$ be the union of all the $2$-handles and $3$-handles in this special handle decomposition of $N$, and let $Z= H''\cup C$. For each $0$-handle $h_i^0$, let $R_i = h_i^0\cap Z$ and $R=\bigcup_i R_i$. 


   For $ 1 \leq i \leq n-1$, let $\cT_i$ be a copy of our new spatial graph tangle in the ball $B_i = h_i^0$ such that the four outgoing ends are identified with endpoints of the core curves of the four $1$-handles that touch $h_i^0$.  Let the exterior of $\cT_i$ in $B_i$ be denoted $Q_i$. 
  
   For each $i$ such that $n \leq i \leq m$, let  $T_i$ be a copy of the true lover's tangle in the ball $B_i = h_i^0$ such that the four outgoing ends are identified with endpoints of the core curves of the four $1$-handles that touch $h_i^0$. In this case, let the exterior of $T_i$ in $B_i$ be denoted $Q'_i$.  We can choose the connections so that the graph $G$, which is the union of all of the tangles in all the 0-handles and the core curves of the 1-handles, is connected.  Let $Q = (\cup_{i =1}^{n-1} Q_i) \cup (\cup_{i = n}^m Q_i')$. Note that a neighborhood of $G$ is a genus $n$ handlebody.

    The exterior $\overline{M-N(G)}$ in $M$ is $Q\cup Z$ with $Q\cap Z=R$. Since $Q$ is a disjoint union of tangle spaces, by work of Myers the tangle space of true lover's tangle with the gluing face has property C$'$, and by Lemma \ref{ELpropertyC'} this also holds for our tangle $J$. It follows that $(Q,R)$ has property C$'$. By Lemma 5.3 of \cite{Myers}, $(Z, R)$ has property B$'$. By Lemma \ref{lemma3.3-Myers}, this means $\overline{M-N(G)}$ is simple. Note that the boundary of the handebody exterior is an incompressible surface in our manifold. Furthermore, the resulting manifold is irreducible and orientable, and thus is Haken. By the work of Thurston it follows that it is tg-hyperbolic.
\end{proof}

The resulting manifold is still compact, so we can repeat the process any number of times, allowing us to remove any sequence of handlebodies and keep the complement hyperbolic. Thus, we obtain the corollary as a general case.

\generalthm*

\section{Volume bounds}\label{volumebounds}


In this section, we consider embeddings of spatial graphs in some elementary 
3-manifolds, and generalize the octahedral decomposition, first detailed by Dylan Thurston in \cite{DThurston} for links in the 3-sphere, to obtain upper bounds on their volumes. The octahedral decomposition of link complements in $S^3$ has already been generalized to links in thickened surfaces, as in \cite{AdamsCalderonMayer}. 

First, we note that the complement of any embedded handlebody $H$ in a 3-manifold $M$ is homeomorphic to the complement of an open regular neighborhood of a deformation retract of $H$, so the theory of exteriors of spatial graphs in 3-manifolds is equivalent to the the theory of closures of complements of handlebodies in 3-manifolds. 

We work with a spatial graph $G$ in $S^3$ or in a thickened surface $F \times I$. We consider spatial graphs with no vertices of degree 1 and with minimum degree of at least 3 so they are not just links. But we do allow more than one component and we allow components that are link components. For spatial graphs in $S^3$, we take $F$ to be a sphere in $S^3$.  For spatial graphs in  $F\times I$, we define projection to $F$ in the usual way. 

If $F$ is a disk, then we recall that to put a hyperbolic metric on a 3-manifold we must cap off all spherical boundaries with 3-balls, so a spatial graph in a 3-ball given by $F \times I$  is equivalent to one in $S^3$. If $F$ has at least one boundary component and is not a disk, then $F\times I$ can be realized as a handlebody. In this case, we still project to $F$. 

For spatial graphs in $S^3$ or a thickened surface, we define crossing number as the minimal number of crossings in a projection to $F$. 


To construct the octahedral decomposition of the complement of a spatial graph, we require each cycle in the graph to have an edge with a crossing. The following proposition says that if this is not the case then the graph is not hyperbolic, so it is reasonable to exclude this case from our decomposition.

\begin{lemma}\label{cyclescrossings}
    Let $M$ be either the 3-sphere $S^3$ or a thickened surface $F \times I$ where $F$ is a compact orientable surface with or without boundary.  If a spatial graph $G$ embedded in $M$ has a projection to $F$ in which a cycle in the graph is involved in no crossings then the exterior of the graph is not hyperbolic.
\end{lemma}
\begin{proof}
    If the cycle without crossings, denoted $L$, is contractible then it bounds an essential disk, precluding hyperbolicity. If $M=S^3$, $L$ is always contractible. If $M$ is not $S^3$ then we may isotope $L$ so that it lies on $\partial M$.
    This isotopy traces out an essential annulus with one boundary on $L$ and the other on $\partial M$ and therefore also precludes hyperbolicity.
\end{proof}

\begin{figure}
    \centering
    \includegraphics[width=0.8\linewidth]{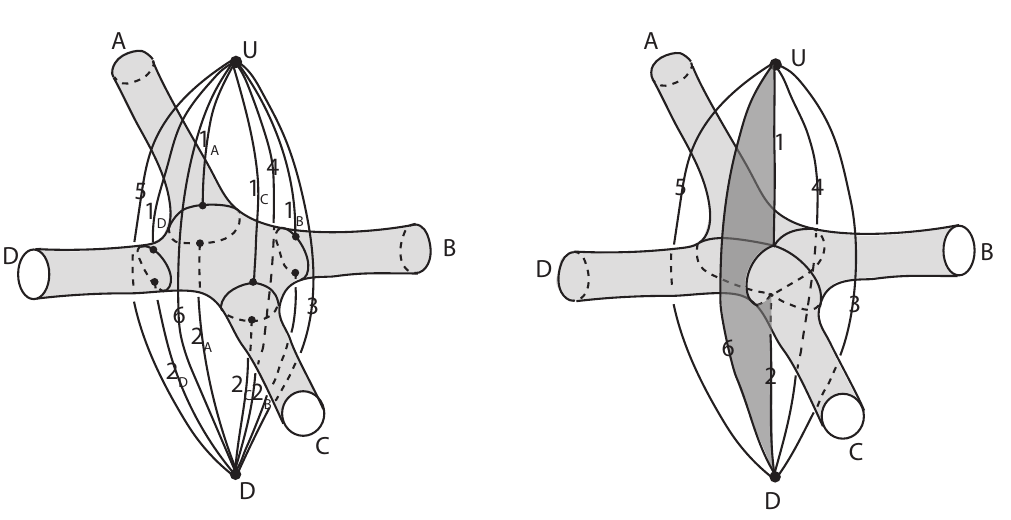}
    \caption{Decomposing manifold around a vertex of the projection of a graph}
    \label{fig:starfruit}
\end{figure}

\begin{construction}[Generalized octahedral decomposition]\label{con:gen-oct-decomp}
Consider a spatial graph $G$  embedded in $S^3$ or a thickened surface $F \times I$ and let $p(G)$ be a minimal crossing projection to $F$. As in the usual octahedral decomposition depicted in Figure \ref{fig:octahedronnew}, we place an octahedron between each crossing of the projection such that one apex is at the overstrand, the other is at the understrand, and two nonadjacent equatorial vertices sit directly below the overstrand with the other two sitting directly above the understrand. We pull the equatorial vertices below the overstrand down to a finite point labeled $D$ below the plane of the projection. It's then clear that the two edges from these equatorial vertices to the understrand will be identified. Similarly, we pull the equatorial vertices above the understrand 
 up to a finite point labeled $U$ above the plane of the projection. This causes the two edges connecting from these equatorial vertices to the overstrand to be identified. The edges of the octahedron at a crossing are identified to edges as in Figure \ref{fig:octahedronnew}. 

\begin{figure}
    \centering
    \includegraphics[width=0.4\linewidth]{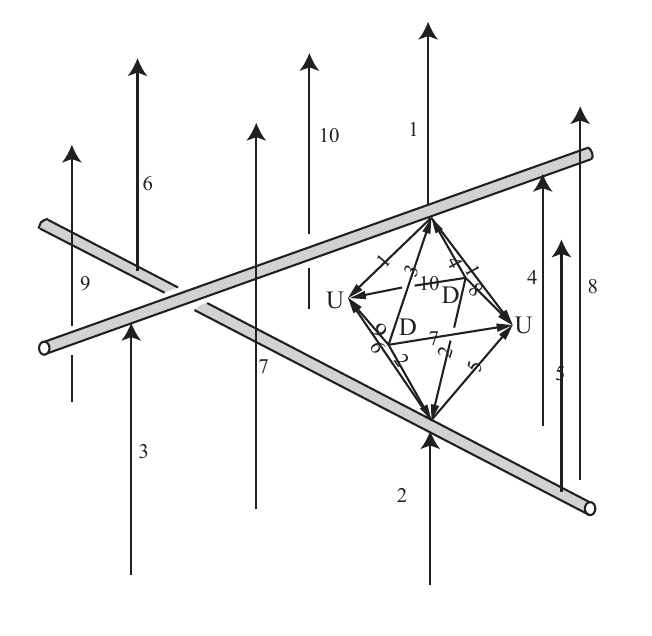}
    \caption{Place one octahedron between each crossing}
    \label{fig:octahedronnew}
\end{figure}

We now need to cover the space surrounding a vertex $v$ of the graph projection. Consider the image on the left in Figure \ref{fig:starfruit}. We see a ball defined by the edges shown around the neighborhood of the vertex. Note that the shaded portion in the figure, which is the neighborhood of $v$ and its outgoing edges, is removed from the ball. Then we isotope each of the four edges labeled $1_A, 1_B, 1_C$, and $1_D$ to a single edge running from the vertex $U$ to the top of the handlebody that is the neighborhood of the graph. Similarly, we isotope the four edges labeled $2_A, 2_B, 2_C$, and $2_D$ to a single edge running from the bottom of the handlebody to the vertex $D$. The ball now becomes an object $Q$ that consists of four fins, one of which is shaded darker in the figure on the right. Each side of a fin corresponds to a triangle that is glued to either the face of an octahedron corresponding to a crossing or to another fin from a vertex that shares an edge with this vertex such that that edge is involved in no crossings. Note that in the case $M \setminus N(G)$ is hyperbolic, Lemma \ref{cyclescrossings} avoids the situation that there is a cycle of fins that glue one to the other and then back to the beginning with no octahedra involved. 

Given a vertex of valency $r$, there will be $r$ such fins. We call this set of fins a {\bf starfruit}.

In the case that the ambient manifold is $S^3$, the vertices $U$ and $D$ remain in the manifold and are finite vertices. 

In the case the manifold is $F \times I$ and the surface $F$ is neither a sphere, disk, torus, nor annulus, the genus of the boundaries that result from $F \times I$ is at least 2. The boundary of the manifold is either two surfaces when $F$ is closed or one surface when $F$ has boundary.  Then, when we take the exterior of the spatial graph, we obtain an additional boundary surface of genus greater than 1 and possibly additional boundary components if $G$ is not connected.  If all boundary components have genus greater than 1, all vertices of both the octahedra and the starfruit are truncated. However, if there are torus boundaries, either occurring when $F$ is an annulus (so there is one torus boundary) or torus (so there are two torus boundaries), or when there are components of $G$ that have no vertices of valency greater than 2, then the corresponding vertices will be ideal.




\end{construction}

The starfruit corresponding to $v$ contributes zero volume to the knot complement. It just changes the gluings on faces of octahedra and the octahedra generate all the volume. Hence, if the spatial graph is hyperbolic, we obtain an upper bound on the volume of its complement, which is the number of crossings times the maximal volume of an octahedron in which every vertex is truncated.



To each crossing, we associate a generalized hyperbolic octahedron, which has vertices that are real, ideal, and hyperideal depending on the genera of $F$ and the connected components of the graph containing the strands which form the crossing. In particular, we remove toroidal boundary components and leave higher genus boundary components. Respectively, these correspond to ideal and hyperideal (also, truncated or ultra-ideal) vertices. Note that hyperideal vertices may be truncated canonically along their truncation planes.

To determine the maximal volume of a generalized octahedron, we use a result of \cite{Belletti2020TheMV} for \textit{proper} generalized polyhedra. The plane along which a vertex of a hyperbolic polyhedron is truncated cuts $\H^3$ into two parts. A \textit{proper} generalized polyhedron is a generalized polyhedron in $\H^3$ for which each truncation plane cuts $\H^3$ such that the side of the plane away from the truncation vertex contains all the real 
vertices of the polyhedron.  Notice that all mildly truncated polyhedra are proper. Also notice that if we consider spatial graphs (possibly with link components) in thickened surfaces, then Construction \ref{con:gen-oct-decomp} yields generalized octahedra with no real vertices. Hence, these octahedra must be proper, and we can use the result below to bound their volumes.

\begin{theorem}[\cite{Belletti2020TheMV}, Theorem 4.2]
    For any planar 3-connected graph $\Gamma$, the supremum of the volumes of all proper generalized hyperbolic polyhedra with 1-skeleton $\Gamma$ is the volume of $\overline{\Gamma}$ the rectification $\Gamma$, that is, the polyhedron in $\R^3$ with 1-skeleton $\Gamma$ all of whose edges are tangent to $\partial\H^3$.
\end{theorem}
\begin{figure}
    \centering
    \includegraphics[width=0.5\linewidth]{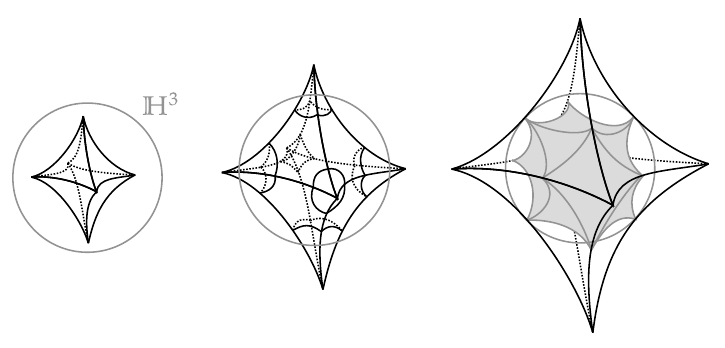}
    \caption{Constructing the rectification of the octahedron}
    \label{fig:oct-rectification}
\end{figure}

It is also shown in \cite{Belletti2020TheMV} that the rectification $\overline{\Gamma}$ of a 3-connected, planar graph $\Gamma$ exists and is unique up to isometry. Moreover, $\overline{\Gamma}$ has a unique associated truncation (and thus a unique volume) which results in the ideal, right-angled hyperbolic polyhedron with 1-skeleton given by the medial graph $M(\Gamma)$. 


\begin{theorem} \label{thm:generalvolbound}
    Let $G$ be a spatial graph embedded in the interior of  a thickened surface $M = F \times I$ where $F$ is a compact orientable surface with or without boundary of Euler characteristic less than 1. Suppose the exterior of  $G$ is tg-hyperbolic. 
    Then, the volume  of the exterior of $G$ is bounded with respect to the crossing number $c(G)$ by 
    $$\mathrm{vol} (\overline{M \setminus N(G)}) < \mathrm{vol}(Q_{cuboct}) \cdot c(G).$$
\end{theorem}

\begin{proof}
     As in the above construction, we have a topological decomposition of the graph exterior into generalized octahedra such that every vertex is either ideal when the corresponding component of the boundary is a torus, or hyper-ideal when the corresponding boundary is of genus 2 or greater.
     
Let $\Gamma_{oct}$ be the 1-skeleton of the octahedron. This is a planar, 3-connected graph satisfying the hypothesis of Theorem 4.2 of \cite{Belletti2020TheMV}, which says in our case that over all proper, generalized octahedra $P_{oct}$ we have
$$\sup_{P_{oct}} \text{vol}(P_{oct}) = \text{vol}(Q_{cuboct}),$$
where $Q_{cuboct}$ is the ideal, right-angled hyperbolic cuboctohedron. Note that in the projective model $\mathbb{H}^3 \subseteq \mathbb{R}^3 \subseteq \mathbb{RP}^3$, a rectification $\overline{\Gamma}$ of a given 3-connected, planar graph $\Gamma$ is a projective polyhedron which has 1-skeleton isomorphic to $\Gamma$ and every edge tangent to the sphere $\partial \mathbb{H}^3$ at infinity. 

As seen in Figure \ref{fig:oct-rectification} we obtain $Q_{cuboct}$ as the truncation of the rectification $\overline{\Gamma_{oct}}$ determined by $\Gamma_{oct}$. Note that a generalized octahedron does not attain this upper bound in volume as it would require the truncation planes to be tangent to one another which cannot occur in a manifold with totally geodesic boundary. 

Given the closure of the complement of $G$ is tg-hyperbolic, in our generalized octahedral decomposition, the crossing number $c(G)$ corresponds to the minimal number of generalized octahedra which can be used to reconstruct the exterior of $G$ under this process. After performing such a decomposition, the volume is then the sum of the volumes of $c(G)$ many generalized octahedra. Note that when realized geometrically, these generalized octahedra may either be singular or negatively oriented, which generates negative volume, but in any case, each has volume upper bounded by $\text{vol}(Q_{cuboct}).$
\end{proof}

Now, we can determine the volume of the right-angled cuboctahedron, which provides the volume bound, as in Theorem \ref{thm:generalvolbound}. This volume has been previously computed, such as in \cite{AdamsCalderonMayer}, but as we are unaware of any sources that discuss the volume in the exact form using elementary methods, we do so here explicitly with symmetries.

\begin{lemma}
    The volume of the ideal, regular cuboctahedron is
\begin{align*}
    \mathrm{vol}(Q_{cuboct}) &= 8 \Lambda(\pi/2 - \theta) + 16 \Lambda(\theta) - 6\Lambda(2\theta) + \Lambda(4\theta)\\
    &\approx 12.04609204009437764726837862923
\end{align*}
where $\theta = \arctan\left({\sqrt 2}\right).$
\end{lemma}

\begin{proof}
    In the Klein ball model of $\mathbb{H}^3$, the ideal, regular polyhedron $Q_{cuboct}$ is realized as a regular cuboctahedron inscribed in the unit sphere. We know that the Poincaré ball model coincides with the Klein model at the unit sphere at infinity, so we take the coordinates for each ideal vertex and continue to calculate the angle with the Poincaré model. 

    We may decompose our cuboctahedron into 13 ideal tetrahedra as depicted in Figure \ref{fig:cuboct2}. The tetrahedra come in three isometry classes. We specify their dihedral angles in terms of  $\theta$, the dihedral angle in the ideal tetrahedron $v_1v_2v_3v_4$ at the edge $v_1v_3$. There are eight tetrahedra isometric to the tetrahedron determined by $v_1v_2v_3v_4$, all with dihedral angles $\pi/2, \pi/2 - \theta, \theta$. Once we remove these tetrahedra (as in the first list of Figure \ref{fig:cuboct2}), we are left with a rectangular parallelepiped. We decompose it into four external tetrahedra, each with three faces on the exterior of the parallelepiped and one in the interior, and one interior tetrahedron, which is shaded in Figure \ref{fig:cuboct2} (second row, left). The four exterior tetrahedra each have dihedral angles $\theta, \theta,  \pi - 2\theta$. The interior tetrahedron has dihedral angles $\pi - 2\theta, \pi - 2\theta, 4\theta - \pi$. 
    \begin{figure}
        \centering
        \includegraphics[width=0.7\linewidth]{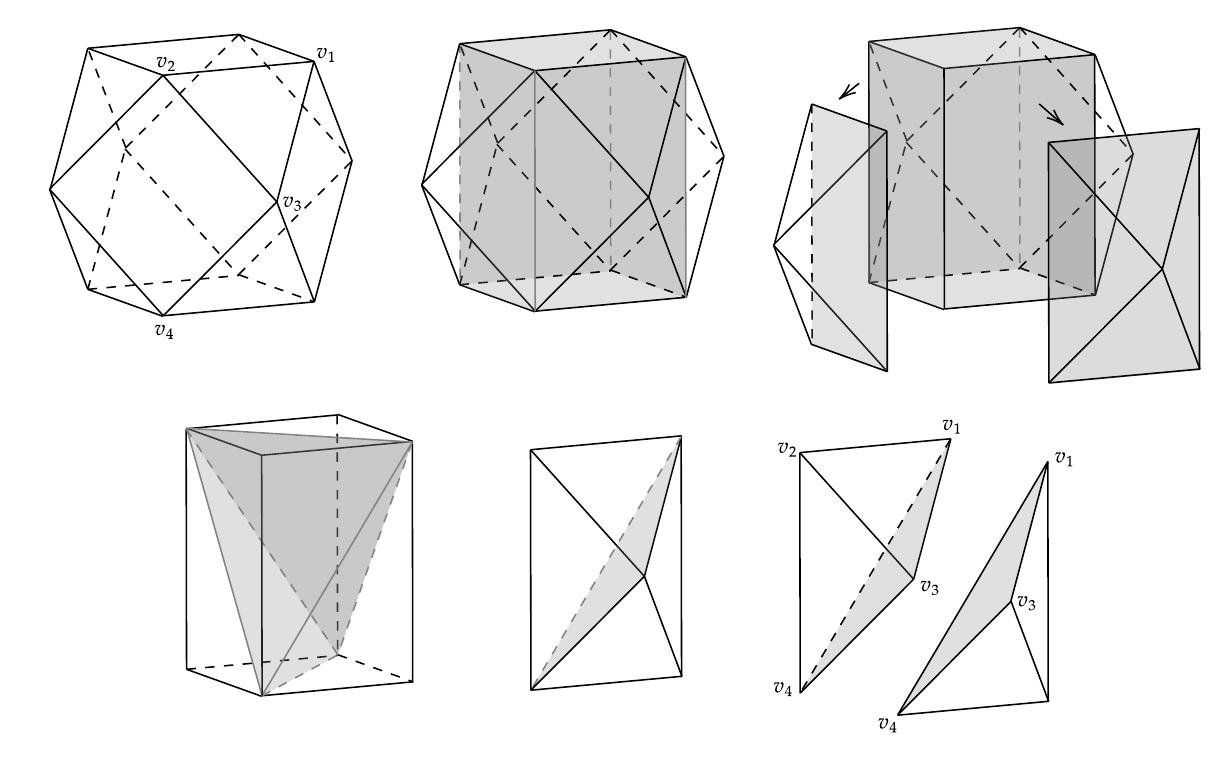}
        \caption{Triangulating the ideal, right angled cuboctahedron $Q_{cuboct}$ in the Klein model}
        \label{fig:cuboct2}
    \end{figure}
    
    In the Poincaré model up to isometry, this tetrahedron can be realized with vertices
    $$v_1 = \left(0, \frac{1}{\sqrt 2}, \frac{1}{\sqrt 2}\right), \quad v_2 = \left(\frac{1}{\sqrt 2}, 0, \frac{1}{\sqrt 2}\right), \quad v_3 = \left(\frac{1}{\sqrt 2}, \frac{1}{\sqrt 2}, 0\right), \quad v_4 = \left(\frac{1}{\sqrt 2}, 0, -\frac{1}{\sqrt 2}\right).$$
    We know that the faces $v_1v_2v_3$ and $v_1v_3v_4$ are on the spheres given by (centered at $c = \left(\frac{1}{\sqrt 2}, \frac{1}{\sqrt 2}, \frac{1}{\sqrt 2}\right)$ and infinity)
    \begin{align*}
        \left(x - \frac{1}{\sqrt 2}\right)^2 + \left(y - \frac{1}{\sqrt 2}\right)^2 + \left(z - \frac{1}{\sqrt 2}\right)^2 &= \frac{1}{2},\\ 
        x - y + z &= 0,
    \end{align*}
    respectively. Their angle of intersection, which by construction is $\theta$, can be found by intersecting these surfaces with some plane, say $x-z = 0$, containing the origin that is orthogonal to both spheres bounding their respective faces. We see that the midpoint $u_1 = \left(\frac{1}{3\sqrt 2}, \frac{\sqrt 2}{3}, \frac{1}{3\sqrt 2}\right)$ of the geodesic edge $v_1v_3$ and $u_2 = \left(\frac{1}{\sqrt 2}, \sqrt 2, \frac{1}{\sqrt 2}\right)$ are the two points on the intersection of all three surfaces. As in Figure \ref{fig:cuboct1}, we can use the Pythagorean theorem to determine the Euclidean length $|u_1u_2| = \frac{2}{\sqrt 3}$ and thus $$\theta = \arctan(\sqrt 2).$$

    \begin{figure}
        \centering
        \includegraphics[width=0.5\linewidth]{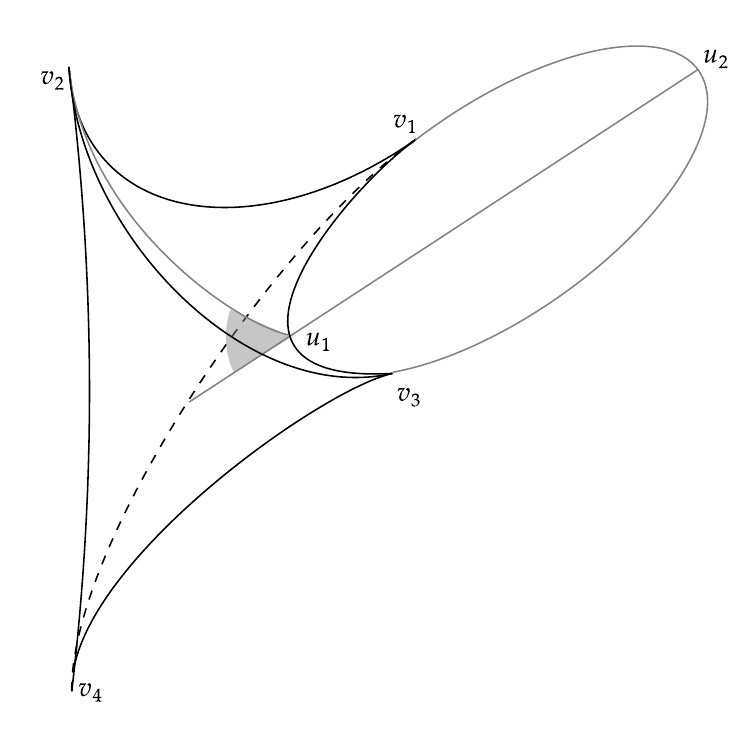}
        \caption{We use the labeled points in the Poincaré model to calculate $\theta$, the shaded angle, which is the dihedral angle.}
        \label{fig:cuboct1}
    \end{figure}
    
    So, we may apply Milnor's formula for the volume of ideal tetrahedra \cite{Milnor150}, where $\Lambda$ denotes the Lobachevsky function:
    \begin{align*}
        \mathrm{vol}(Q_{cuboct}) &= 8\left(\Lambda(\theta) + \Lambda(\pi/2 - \theta) \right) + 4 \left(2\Lambda(\theta) + \Lambda(\pi - 2\theta)\right) + \left(2\Lambda(\pi - 2\theta) + \Lambda(4\theta - \pi)\right)\\
        &= 8 \Lambda(\pi/2 - \theta) + 16 \Lambda(\theta) + 6 \Lambda(\pi - 2\theta) + \Lambda(4\theta - \pi)\\
        &= 8 \Lambda(\pi/2 - \theta) + 16 \Lambda(\theta) - 6\Lambda(2\theta) + \Lambda(4\theta),
    \end{align*}
    and we have the volume in its exact form.
\end{proof}

By the work of Igor Rivin, the cuboctahedron $Q_{cuboct}$ is unique up to isometry and, being maximally symmetric, has maximal volume among all ideal, convex cuboctahedra \cite{RivinMain} [Theorems 14.1, 14.3].



\begin{theorem}
Let $G$ be a spatial graph embedded in $S^3$ with projection to a sphere $F$ in $S^3$ with $c(G)$ crossings.  Suppose the exterior of  $G$ is tg-hyperbolic. 
    Then, the volume  of the exterior of $G$ is bounded with respect to the crossing number $c(G)$ by 
    $$\mathrm{vol} (\overline{S^3 \setminus N(G)}) < \mathrm{vol}(B^{\mathrm{trunc}}_4) \cdot c(G),$$
    where $B^{\mathrm{trunc}}_4$ is the 4-bipyramid of maximal volume described in \cite{AdamsCalderonMayer}, with volume approximately 5.074708.
\end{theorem}

\begin{proof}
    We work in the Klein model. In the generalized decomposition for the complement in $S^3,$ we obtain a generalized octahedron $P$ with two opposite vertices $u_1, u_2$ that are either ideal or hyperideal, determined by the genera of the handlebody components forming the crossing, and four real equatorial vertices $v_1, v_2, v_3, v_4$, connected by edges cyclically in this order. 

    Let $v_1'$ and $v_3'$ be the two ideal endpoints of the unique geodesic through $v_1$ and $v_3$. Let $v_2'$ and $v_4'$ be the two ideal endpoints of the unique geodesic through $v_2$ and $v_4$. And let $P'$ be the generalized octahedron with vertices $u_1, u_2, v_1', v_2', v_3'$ and $v_4'$.
    Then $P \subset P'$.

    
    Since every such generalized octahedron $P$ is contained in a generalized octahedron (i.e. 4-bipyramid) $P'$ with four equatorial ideal vertices, we may bound $\mathrm{vol}(P)$ by the maximum volume of such an octahedron $P'$. In \cite{AdamsCalderonMayer}, this volume is found to be $$\mathrm{vol}(B^{\mathrm{trunc}}_4)  \approx 5.07470803204826812510601277,$$ realized by $B^{\mathrm{trunc}}_4$, the generalized octahedron with dihedral angle $2\pi/3$ on the four equatorial edges and dihedral angle $\pi/3$ on each of the other eight edges. As in the proof of the previous theorem, when realized geometrically, these generalized octahedra can be either singular or negatively oriented, in which case they contribute negative volume. But, in the case of each generalized octahedron, the volume contribution is bounded strictly from above by the upper bound, yielding the conclusion of the theorem. 
\end{proof}

In fact, results such as these can yield volume bounds for handlebody complements in much more general 3-manifolds. We utilize the following two theorems. The first is well known \cite[Thm 1.10.15]{Klingenberg}.

\begin{theorem}\label{fixedsurfacelemma}
    Let $M$ be a Riemannian manifold and let $f:M \rightarrow M$ be a nontrivial isometry.
    Then, the subset Fix$(f)$ is a union of embedded totally geodesic submanifolds.
\end{theorem}

  The second is due to Agol, Storm, and Thurston in  \cite[Thm. 9.1]{Agol-Storm-Thurston}, as  generalized by Calegari, Freedman and Walker in \cite[Thm. 5.5]{CFW}.
  
\begin{theorem}\label{cuttingandpastingthm}
    Let $\overline{M}$ be a compact manifold with interior $M$ a hyperbolic 3-manifold of finite-volume.
    Let $\Sigma$ be a properly embedded two-sided totally geodesic surface, let $\psi:\Sigma \rightarrow \Sigma$ be a diffeomorphism, and let $M'$ be the manifold formed by surgering along $\Sigma$ and gluing the resulting pieces together via $\psi$.
        Then, $M'$ is a hyperbolic 3-manifold of finite volume, satisfying
        \[
            \vol(M') \geq \vol(M).
        \]
        Equality is attained if and only if $\Sigma$ is totally geodesic in $M'$.
\end{theorem}

In our case, we will apply these as follows. Let $M$ be a compact orientable 3-manifold and let $F$ be a properly embedded incompressible two-sided compact orientable surface in $M$.  Let $D(M \backslash \backslash F)$ be the double of the manifold $M$ obtained by removing an open neighborhood of $F$, and then doubling $M \backslash \backslash F$ across the two copies of $F$ in its boundary. Note that by Theorem \ref{fixedsurfacelemma}, the two copies of $F$ in the double are totally geodesic, as there is a reflection of $D(M \backslash \backslash F)$ for which they are the fixed point set. Let $D(\overline{F \times I \setminus N(G)})$ denote the double of $\overline{F \times I \setminus N(G)}$ over the two copies of $F$ on its boundary. For it also, the two relevant copies of $F$ in $D(\overline{F \times I \setminus N(G)})$ are also totally geodesic. Note that when $\overline{F \times I \setminus N(G)}$ is tg-hyperbolic, so does $D(\overline{F \times I \setminus N(G)})$, and has volume exactly twice the volume of $\overline{F \times I \setminus N(G)}$.

\begin{theorem}\label{genmanthm} Let $M$ be a compact orientable 3-manifold with properly embedded incompressible compact orientable two-sided surface $F$.  Let $G$ be a spatial graph in the interior of $F \times I \subset M$. If $D(M \backslash \backslash F)$ and $\overline{F \times I \setminus N(G)}$ are both tg-hyperbolic, then $\overline{M \setminus N(G)}$ is tg-hyperbolic and $$\mathrm{vol}(\overline{M \setminus N(G)}) \geq \frac{1}{2} \mathrm{vol} (D(M \backslash \backslash F)) + \mathrm{vol}(\overline{F \times I \setminus N(G))}.$$
 \end{theorem}

 Note that the original surface $F$ can have boundary or not and it need not be totally geodesic in $M$. But if it is, then $$\frac{1}{2} \mathrm{vol} (D(M \backslash \backslash F)) = \mathrm{vol} (M).$$

\begin{proof} We start with the disconnected manifold $W = D(M \backslash \backslash F) \cup D(\overline{F \times I \setminus N(G))}$. Denote the two totally geodesic copies of $F$ in $D(M \backslash \backslash F)$ by $F_0^A$ and $F_1^A$. 
 Denote the two totally geodesic copies of $F$ in the double of $\overline{F \times I \setminus N(G)}$ by $F_0^B$ And $F_1^B$.   Cutting  $D(M \backslash \backslash F)$ open along $F_0^A$ And $F_1^A$ yields two copies of $M \backslash \backslash F$, the first with copies $F_{0, 0}^A$ and $F_{1, 0}^A$ of $F$ and the second with copies $F_{0, 1}^A$ and $F_{1, 1}^A$ of $F$. Cutting $D(\overline{F \times I \setminus N(G))}$ open along $F_0^B$ And $F_1^B$ yields two copies of $\overline{F \times I \setminus N(G)}$, the first with copies $F_{0, 0}^B$ and $F_{1, 0}^B$ of $F$ and the second with copies $F_{0, 1}^B$ and $F_{1, 1}^B$ of $F$. These cuts have the effect of cutting $W$ open along four totally geodesic copies of $F$ and yield the eight listed copies in the resulting manifold. Now, we glue  $F_{i, j}^A$ to $F_{i, j}^B$,  for all $ i = 0, 1$ and $j = 0, 1$ in each case using the same diffeomorphism that came from how they were glued before the cutting. This creates two copies of $\overline{M \setminus N(G)}$. By Theorem \ref{cuttingandpastingthm}, the volume of the resulting two copies of $\overline{M \setminus N(G)}$ is at least as large as the volume of $W$. Dividing by 2 yields the desired inequality.  
     \end{proof}

\begin{figure}
        \centering
        \includegraphics[width=0.5\linewidth]{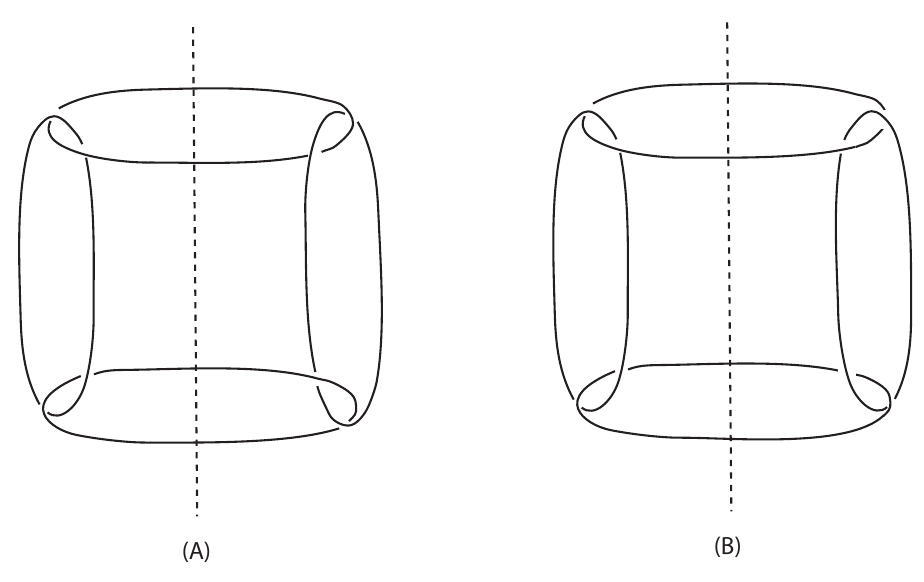}
        \caption{A hyperbolic link complement $M$ containing an incompressible surface $F$ such that $D(M \backslash \backslash F)$ is not hyperbolic}
        \label{fig:doublecounter}
    \end{figure}

 Note that if $M$ is tg-hyperbolic, this needn't imply that $D(M \backslash \backslash F)$ is tg-hyperbolic. As a simple example, take the exterior of an alternating 4-chain link as in Figure \ref{fig:doublecounter}(A), which is tg-hyperbolic. Then let $F$ be the 4-punctured sphere corresponding to the dotted line. Then $D(M \backslash \backslash F)$ is the exterior of the link in Figure \ref{fig:doublecounter}(B), which contains essential annuli and is therefore not tg-hyperbolic.
 
 However, when $D(M \backslash \backslash F)$ is tg-hyperbolic, knowledge of the volume of $\overline{F \times I \setminus N(G)}$ allows us to obtain lower bounds on the volume of $\overline{M \setminus N(G)}$. 

\bibliographystyle{plain}

\bibliography{ref}

\end{document}